\documentclass{amsart}
\usepackage{amsmath,amssymb,amscd,amsfonts}

\newtheorem{theorem}{Theorem}
\newcommand{\bt}{\begin{theorem}}
\newcommand{\et}{\end{theorem}}
\newtheorem{lemma}{Lemma}
\newcommand{\bl}{\begin{lemma}}
\newcommand{\el}{\end{lemma}}
\newtheorem{corollary}{Corollary}
\newcommand{\bc}{\begin{corollary}}
\newcommand{\ec}{\end{corollary}}
\newcommand{\beq}{\begin{equation}}
\newcommand{\eeq}{\end{equation}}
\newcommand{\benum}{\begin{enumerate}}
\newcommand{\eenum}{\end{enumerate}}
\newcommand{\N}{\ensuremath{ \mathbf N }}
\newcommand{\Z}{\ensuremath{\mathbf Z}}

\newcommand{\mca}{\ensuremath{ \mathcal A}}
\newcommand{\mcb}{\ensuremath{ \mathcal B}}

\newcommand{\mcg}{\ensuremath{ \mathcal G}}

\DeclareMathOperator{\qqand}{\qquad\text{and}\qquad}

\title{Limits and decomposition of 
de Bruijn's additive systems}

\author{Melvyn B. Nathanson}
\address{Lehman College (CUNY),Bronx, New York 10468}
\email{melvyn.nathanson@lehman.cuny.edu}

\thanks{Supported in part by a grant from the PSC-CUNY Research Award Program.}

\date{\today}

\begin{document}

\begin{abstract}
An \emph{additive system} for the nonnegative integers is a 
family $(A_i)_{i\in I}$ of sets of nonnegative integers 
with $0 \in A_i$ for all $i \in I$ such that every nonnegative integer can be written 
uniquely in the form $\sum_{i\in I} a_i$ with $a_i \in A_i$ for all $i$ 
and $a_i \neq 0$ for only finitely many $i$.  
In 1956, de Bruijn proved that every additive system is constructed
 from  an infinite sequence $(g_i )_{i \in \N}$  of integers with $g_i \geq 2$ for all $i$, 
 or is a contraction of such a system.  
This paper discusses limits and the stability of additive systems, 
and also describes the ``uncontractable'' or 
``indecomposable'' additive systems.
\end{abstract}

\maketitle

\section{Additive systems and de Bruijn's theorem}

Let $\N_0$ and \N\ denote the sets of nonnegative integers and positive
integers, respectively.  
For real numbers $a$ and $b$, we define the intervals of integers
 $[a,b) = \{x\in \Z : a \leq x < b\}$
 and  $[a,b] = \{x\in \Z : a \leq x \leq b\}$.   

Let $I$ be a nonempty finite or infinite set, 
and let $\mca = (A_i)_{i\in I}$ be a family of  sets of integers with $0 \in A_i$ and 
$|A_i| \geq 2$ for all $i \in I$.   
Each set $A_i$ can be finite or infinite.  
The \emph{sumset} $S = \sum_{i\in I} A_i$ is the set of all integers $n$ 
that can be represented in the form $n = \sum_{i\in I} a_i$, 
where $a_i \in A_i$ for all $i \in I$ and $a_i  \neq 0$ for only finitely many $i \in I$.  
If every element of  $S$ has a \emph{unique} 
representation in the form $n = \sum_{i\in I} a_i$, 
then we call \mca\ a  \emph{unique representation system for $S$}, 
\index{unique representation system}
and we write 
$
S =  \bigoplus_{i\in I} A_i.
$

In a unique representation system  \mca\  for $S$, 
we have  $A_i \cap A_j = \{ 0 \}$ for all $i \neq j$. 
The condition $|A_i| \geq 2$ for all $i\in I$ implies that  
$A_i = S$ for some $i\in I$ if and only if $|I|=1$.  
Moreover, if $I^{\flat} \subseteq I$ and $S =  \sum_{i\in I^{\flat}} A_i$, 
then $S =  \bigoplus_{i\in I^{\flat}} A_i$ and $I = I^{\flat}$.

The family $\mca = (A_i)_{i\in I}$ is  an \emph{additive system}
\index{additive system} if \mca\  is a unique representation system 
for the set of nonnegative integers, that is, 
if $\N_0 =  \bigoplus_{i\in I} A_i$.  
The following lemma follows immediately from the definition of an additive system.

\bl     \label{deBruijn:lemma:contraction}
Let $\mcb = (B_j )_{j\in J} $  be an additive system.  
If $\{J_i\}_{i\in I}$ is a partition of $J$ into pairwise disjoint nonempty sets, and  if 
\[
A_i= \sum_{j \in J_i} B_j
\]
then $\mca = (A_i )_{i\in I}$ is an additive system.
\el

The additive system \mca\ obtained from the additive system \mcb\ 
by the partition procedure described in 
Lemma~\ref{deBruijn:lemma:contraction} 
is called a \emph{contraction} of \mcb.  
(In~\cite{debr56},  de Bruijn called \mca\ a  \emph{degeneration}  
\index{degenerate additive system} of \mcb.) 
If $I=J$ and if $\sigma$ is a permutation of $J$ such that $J_i =\{\sigma(i) \}$
for all $i \in J$, then \mca\ and \mcb\ contain exactly the same sets.   
Thus, every additive system is a contraction of itself.  
An additive system \mca\ is a \emph{proper contraction} of \mcb\
if at least one set $A_i \in \mca$ is the sum of at least two sets in \mcb.

Let $X$ be a set of integers and let $g$ be an integer.  
The \emph{dilation} of $X$ by $g$ is the set $g\ast X = \{g x : x \in X \}$.

\bl     \label{deBruijn:lemma:dilation}
Let $\mcb = (B_j )_{j\in J} $  be an additive system.  
and let $I = \{i_0\} \cup J$, where $i_0 \notin J$.
If
\[
A_{i_0} = [0,g)
\]
and
\[
A_j = g\ast B_j \qquad\text{for all $j \in J$}
\]
then $\mca = (A_i )_{i\in I}$ is an additive system.
\el

The additive system \mca\ obtained from the additive system \mcb\ 
by the procedure described in 
Lemma~\ref{deBruijn:lemma:dilation} 
is called the \emph{dilation} of \mcb\ by $g$.

There are  certain additive systems that de Bruijn called
\index{British number system}   \emph{British number systems}.
A British number system 
is an additive system constructed from an infinite 
sequence of integers according to the algorithm in 
Theorem~\ref{deBruijn:theorem:BNS} below.
de Bruijn~\cite{debr56} proved  that British number systems 
are essentially the only additive systems.  

\bt    \label{deBruijn:theorem:BNS}
Let $(g_i )_{i \in \N}$ be an infinite sequence of integers 
such that $g_i \geq 2$ for all $i \geq 1$.  
Let $G_0 = 1$ and, for $i \in \N$, let  $G_i = \prod_{j=1}^i g_j$ and 
\[
A_i  = \{ 0, G_{i-1}, 2 G_{i-1},\ldots, (g_i-1)G_{i-1} \} = G_{i-1} \ast [0,g_i).
\]
Then $\mca = (A_i )_{i \in \N}$  is an additive system.  
\et

\bt    \label{deBruijn:theorem:BNS2}
Every additive system is a British number system 
or a proper contraction of a British number system.
\et

The proof of Theorem~\ref{deBruijn:theorem:BNS2} depends on the following fundamental lemma.

\bl     \label{deBruijn:lemma:deBruijn}
Let  $\mca = (A_i)_{i\in I}$ be an additive system with $|I| \geq 2$,
and let $i_1$ be the unique element of $I$ such that $1 \in I_{i_1}$.
There exist an integer $g \geq 2$   
and a family of sets $\mcb = (B_i)_{i\in I}$ such that 
\[
A_{i_1} = [0, g)  \oplus g\ast B_{i_1}
\]
and, for all $i \in I \setminus \{i_1\}$,
\[
A_i = g\ast B_i.
\]
If $B_{i_1} = \{ 0\}$, then $\mcb = (B_i)_{i\in I\setminus \{i_1\}}$ is an additive system, 
and \mca\ is the dilation of the additive system \mcb\ by the integer $g$. 
If $B_{i_1} \neq \{ 0\}$, then $\mcb = (B_i)_{i\in I}$ is an additive system 
and \mca\ is a contraction of the additive system \mcb\ dilated by $g$.
\el

For proofs of Lemmas~\ref{deBruijn:lemma:contraction}, 
\ref{deBruijn:lemma:dilation}, and~\ref{deBruijn:lemma:deBruijn} 
and Theorems~\ref{deBruijn:theorem:BNS} 
and~\ref{deBruijn:theorem:BNS2}, see Nathanson~\cite{nath13}.

This paper gives a refinement of de Bruijn's theorem.
Every additive system is a contraction of a British number system, 
but even a British number system can be a proper contraction of another 
British number system.  An additive system that is not a proper contraction of another 
number system will be called \emph{indecomposable}.  
In Section~\ref{deBruijn:section:decompose} we describe all 
indecomposable British number systems.
Unsurprisingly, there is a one-to-one correspondence between indecomposable 
British number systems and infinite sequences of prime numbers.  

In Section~\ref{deBruijn:section:limits} we define the limit of a sequence of 
additive systems and discuss the stability of British number systems.

Maltenfort~\cite{malt17} and Munagi~\cite{muna05} have also 
studied de Bruijn's additive systems.

\section{Decomposable and indecomposable sets} 

The set $A$ of  integers is a \emph{proper sumset} \index{sumset} 
 \index{sumset!proper}    \index{proper sumset} 
if there exist sets $B$ and $C$ of integers such that $|B| \geq 2$, 
$|C| \geq 2$, and $A = B + C$.   
For example, if $u$ and $v$ are integers and $v-u \geq 3$, 
then the interval $[u, v)$ is a proper sumset:
\[
[u, v) = [0,i) + [u, v + 1 -i)
\]
for every $i \in [2,v-u)$.

The set $A$ of integers is \emph{decomposable} 
\index{decomposable set} 
if there exist sets $B$ and $C$ such that 
$(B,C)$ is a unique representation system for $A$, 
that is, if $|B| \geq 2$, $|C| \geq 2$, and $A=B \oplus C$. 
A decomposition $A = B \oplus C$ is also called a \emph{tiling}\index{tiling} of $A$ by $B$.  
For example, 
\[
[0,12) = \{0,3\} \oplus \{0,1,2,6,7,8\}.
\]
If $A = B \oplus C$ is a decomposition, then $|A| = |B| \ |C|$ 
and so the integer $|A|$ is composite.

Let $n \geq 2$ and consider the interval of integers $A = [0,n)$.
A \emph{proper divisor}\index{proper divisor} of $n$ 
is a divisor $d$ of $n$ such that $1 < d < n$.  
Associated to every proper divisor $d$ of $n$ is the decomposition  
\beq              \label{deBruijn:d-decomp}
[0,n) = [0, d) \oplus d\ast [0, n/d).
\eeq
This is simply the division algorithm for integers.  
The number of decompositions of type~\eqref{deBruijn:d-decomp} 
is the number of proper divisors $d$ of $n$.  
There is exactly one such decomposition if and only if the integer $n$ has a unique 
proper divisor if and only if $n$ is the square of a prime number.

\bl 
Let $n \geq 2$.  
The interval $[0,n)$ is indecomposable if and only if $n$ is prime.
\el

\begin{proof}
If $n$ is prime then $[0, n)$ is indecomposable, 
and if $n$ is composite, then $[0, n)$ is decomposable.  
\end{proof}

If $A=B \oplus C$ and $g$ is a nonzero integer, 
then $g \ast A= g \ast B \oplus g \ast C$, 
and so every dilation of a decomposable set is decomposable.  

The \emph{translate}\index{translate} of the set $A$ by an integer $t$ is the set
\[
A+t = \{a+t:a\in A \}.
\]
Let $t_1,t_2 \in \Z$ with $t = t_1+t_2$.  If $A = B+C$, then 
\[
A+t = (B+t_1) + (C+t_2).
\]
In particular, $A+t = (B+t) + C$.  
If $A=B \oplus C$, then $A + t=  (B + t)  \oplus C$, 
and so every translate of a decomposable set is decomposable.  
Similarly, if $A=B \oplus C$, then $A =  (B - t)  \oplus (C+t)$ 
for every integer $t$.  

Let $A$ be a set of nonnegative integers with $0 \in A$, and let $B$ and $C$ be 
sets of integers with $A = B \oplus C$.  Let $t = \min(B)$.  
Defining $B' = B-t$ and $C' = C+t$, we obtain 
$A = B' \oplus C'$.  Because $\min( B') = 0$, we obtain 
\[
0 = \min(A) = \min(B') + \min(C') = 0 + \min(C') = \min(C')
\]
and so $B'$ and $C'$ are sets of nonnegative integers with $0 \in B' \cap C'$.  

Not every set with a composite number of elements is decomposable. 
For example, the $n$-element set $\{0,1,2,2^2,\ldots, 2^{n-2}\}$ 
is indecomposable for every $n \geq 2$.  This is a special case of the following result.

\bl                             \label{deBruijn:lemma:2Classes}
Let $m \geq 2$.  
Let $A$ be a set of integers that contains integers $a_0$ and $a_1$ 
such that  $a_0 \not\equiv a_1 \pmod{m}$, 
and $a \equiv a_0 \pmod{m}$ for all $a \in A\setminus \{a_1\}$.  The set $A$ is indecomposable.   
\el

\begin{proof}
The distinct congruence classes $a_0 \pmod{m}$ and $a_1\pmod{m}$ contain elements of $A$.   
Let $B$ and $C$ be sets of integers such that $A = B +C$ with $|B|, |C| \geq 2$.
If $B$ is contained in the congruence class $r \pmod{m}$ and 
$C$ is contained in the congruence class $s \pmod{m}$, then $B+C$ 
is contained in the congruence class $r+s \pmod{m}$, and so $A \neq B+C$ 
(because $A$ intersects two congruence classes).
Therefore, at least one of the sets $B$ and $C$ must contain 
elements from distinct congruence classes modulo $m$.  
Let $b_1,b_2 \in B$ with $b_1 \not\equiv  b_2 \pmod{m}$, and 
let $c_1,c_2 \in C$ with $c_1\neq c_2$.  
We have $b_i+c_1 \in B+C$ for $i = 1,2$ and $b_1+c_1 \not\equiv  b_2 + c_1 \pmod{m}$.  
Because $A$ intersects only two congruence classes modulo $m$, and because the 
intersection with the congruence class $a_1 \pmod{m}$ contains only the integer $a_1$, 
we must have $b_i+c_1 = a_1$ for some $i \in \{1,2\}$.  

Similarly, $b_j+c_2  \in B+C$ for $j = 1,2$ 
with $b_1+c_2 \not\equiv  b_2 + c_2 \pmod{m}$, 
and so  $b_j+c_2 = a_1$ for some $j \in \{1,2\}$.  
The equation $b_i+c_1 = b_j+c_2$ implies that $A \neq B \oplus C$.
This completes the proof.  
\end{proof}

The following examples show that, in Lemma~\ref{deBruijn:lemma:2Classes}, 
the condition that the set $A$ contains 
exactly one element of the congruence class $a_1 \pmod{m}$ is necessary.  

Let $m \geq 2$, and  let $R \subseteq [0,m)$ with $|R| \geq 2$.
For every set $J$ of integers with $|J| \geq 2$, we have 
\[
A  =  \{jm + r : j\in J \text{ and } r \in R \} = B \oplus C 
\]
where 
\[
B  =  \{jm:j\in J\} \qqand C = R. 
\]

Let $k$, $\ell$, and $m$ be integers with $k \geq 2$, $\ell \geq 2$, and $m \geq 2$, 
and let $u$ and $v$ be integers 
such that $u \not\equiv v \pmod{m}$.  
Consider the set   
\[
A = \{im+u: i \in [0,\ell ) \} \cup \{ j m  + v: j \in [0, k \ell)   \}.
\]
The sets  
\[
B  = \{ u \} \cup \{ q \ell m  + v: q \in [0, k )   \}
\]
and
\[
C = \{im : i \in  [0, \ell ) \}
\]
satisfy $|B| = 1 + k\ell \geq 2$, $|C| = \ell \geq 2$ and 
\[
A = B \oplus C.  
\]

\section{Decomposition  of additive systems}   \label{deBruijn:section:decompose}
Contraction and dilation are two methods to construct new additive systems from old ones.  
Decomposition is a third  method to produce new additive systems.

An additive system $\mca = (A_i)_{i\in I}$ 
is called \emph{decomposable} \index{decomposable system} 
 \index{additive system!decomposable} 
if  the set $A_{i_0}$ is decomposable for some $i_0 \in I$, 
and  \emph{indecomposable} \index{indecomposable system} 
 \index{additive system!indecomposable} 
if $A_i$ is indecomposable for all $i \in I$.  
Equivalently, an indecomposable additive system is an additive system 
that is not a proper contraction of another additive system.

\bt
Let $\mca = (A_i)_{i\in I}$  be a decomposable additive system,  
and let $A_{i_0}$ be a decomposable set in \mca. 
Choose sets $B$ and $C$ of nonnegative integers such that $0 \in B\cap C$, 
$|B| \geq 2$, $|C| \geq 2$, and $A_{i_0} = B \oplus C$.    
Let
\[
I' = \{ j_1,j_2\} \cup I\setminus \{ i_0\}.
\]
The family of sets  
$\mca' = (A'_i)_{i\in I'}$ defined by 
\[
A'_i = 
\begin{cases}
A_i & \text{ if } i \in  I\setminus \{ i_0\} \\
B & \text{ if } i=j_1 \\
C & \text{ if } i= j_2
\end{cases}
\]
is an additive system.
\et

\begin{proof}
This follows immediately from the definitions of additive system and indecomposable set.  
\end{proof}

We call $\mca'$ a \emph{decomposition} 
of the additive system \mca.

\bl        \label{deBruijn:lemma:AX=B}
Let $a$ and $b$ be positive integers, and let $X$ be a set of integers.
Then
\beq                                       \label{deBruijn:AX=B}
[0,ab) = [0,a) \oplus X 
\eeq
if and only if 
\[
X = a\ast [0,b).
\]
\el

\begin{proof}
The division algorithm implies that $ [0,ab) = [0,a) \oplus a\ast [0,b)$, 
and so  $X = a\ast [0,b)$ is a solution of the additive set equation~\eqref{deBruijn:AX=B}.

Conversely, let $X$ be any solution of~\eqref{deBruijn:AX=B}. 
Let $I = \{1,2,3\}$ and let $A_1 = [0,a)$, $A_2 = X$, and $A_3 =  ab\ast N_0$.
By the division algorithm, $\mca = (A_i)_{i\in I}$ 
is an additive system.  Applying Lemma~\ref{deBruijn:lemma:deBruijn} 
to \mca, we obtain an integer $g \geq 2$ and sets $B_1$, $B_2$, 
and $B_3$ such that 
\begin{align*}
 [0,a) & = [0,g) \oplus g \ast B_1 \\
X & = g\ast B_2 \\
 ab\ast \N_0 & = g \ast B_3.
\end{align*}
It follows that $g=a$,  $B_1 = \{0\}$, $B_3 = b\ast \N_0$, 
and 
\[
\N_0 = B_2 \oplus B_3 = B_2 \oplus b \ast \N_0.
\]
This implies that $B_2 = [0,b)$ and $X = a\ast [0,b)$.   
\end{proof}

There is also a nice polynomial proof of Lemma~\ref{deBruijn:lemma:AX=B}.
Let
\begin{align*}
f(t) & = \sum_{i\in [0,ab)} t^i \\
g(t) &  = \sum_{j \in [0,a)} t^j  \\
h(t) &  = \sum_{k\in [0,b)} t^{ak} \\
h_X(t) & = \sum_{x\in X} t^x.
\end{align*}
The set equation $[0,ab) = [0,a) \oplus  a\ast [0,b)$ implies that 
\[
f(t) = g(t) h(t).  
\]
If $[0,ab) = [0,a) \oplus X$, then 
\[
f(t) = g(t) h_X(t)
\]
and so 
\[
g(t)( h(t) - h_X(t)) = 0.
\]
Because $g(t) \neq 0$, it follows that $h(t) = h_X(t)$ or, equivalently, 
$a\ast [0,b) = X$.

By Theorem~\ref{deBruijn:theorem:BNS2}, 
every additive system is a British number system 
or a proper contraction of a British number system.  
However,  a British number system can also be a proper contraction of another 
British number system. 
Consider, for example, the British number systems $\mca_2$ and $\mca_4$ 
generated by the sequences $(2)_{i \in \N}$ and $(4)_{i \in \N}$, respectively:  
\begin{align*}
\mca_2 & = (\{ 0, 2^{i-1} \} )_{i \in \N} = ( 2^{i-1}\ast [0,2) )_{i \in \N} \\
& = \left( \{ 0, 1 \}, \{ 0, 2 \},  \{ 0, 4 \}, \{ 0, 8 \}, \ldots \right) 
\end{align*}
and
\begin{align*}
\mca_4 & =  ( \{ 0, 4^{i-1}, 2\cdot 4^{i-1}, 3\cdot 4^{i-1} \} )_{i \in \N}
= ( 4^{i-1}\ast [0,4) )_{i \in \N} \\
& = \left( \{ 0, 1, 2, 3 \}, \{ 0, 4,8,12 \},  \{ 0, 16, 32, 48 \}, \{ 0, 64,128, 192, 256 \}, \ldots \right). 
\end{align*}
Because  
\[
4^{i-1}\ast [0,4) = \{ 0, 2^{2i-2} \} + \{ 0, 2^{2i-1} \} = 2^{2i-2} \ast [0,2) + 2^{2i-1} \ast [0,2)
\]
we see that $\mca_4$ is a contraction of $\mca_2$.

de Bruijn~\cite{debr56} asserted the following necessary and sufficient 
condition for one British number system to be a contraction 
of another British number system.  

\bt
Let $\mcb = (B_j)_{j \in \N}$ be the British number system constructed 
from the integer sequence $(h_j)_{j \in \N}$, 
and let $\mca = (A_i)_{i \in \N}$ be the contraction of \mcb\ 
constructed from a partition $(J_i)_{i \in \N}$ of \N\ into nonempty finite sets.  
Then \mca\ is a British number system if and only if $J_i$ 
is a finite interval of integers for all $i \in \N$. 
\et

\begin{proof}
Let $(J_i)_{i \in \N}$ be a partition of \N\ into nonempty finite intervals of integers.
After re-indexing, there is a strictly increasing sequence 
$(u_i)_{i \in \N_0}$ of  integers with $u_0 = 0$ such that 
$J_i = [ u_{i-1}+1 , u_i]$ for all $i \in \N$.

If $\mcb = (B_j)_{j \in \N}$ is the British number system constructed 
from the integer sequence $(h_j)_{j \in \N}$, 
then $B_j = H_{j-1} \ast [0,h_j)$, where $H_0=1$ and $H_j = \prod_{k=1}^j h_k$.    
Let $G_0 = 1$.   For $i \in \N$ we define
\[
g_i = \frac{H_{u_i}}{H_{u_{i-1}}}
\]
and 
\[
G_i = \prod_{j=1}^i g_j = \prod_{j=1}^i \frac{H_{u_j}}{H_{u_{j-1}}} = H_{u_i}.
\]
We have
\begin{align*}
A_i & = \bigoplus_{j \in J_i} B_j 
 = \bigoplus_{j = u_{i-1}+1}^{u_i}  H_{j-1}\ast [0,h_j) \\
& = H_{u_{i-1}} \ast \bigoplus_{j = u_{i-1}+1}^{u_i}  \frac{H_{j-1} }{ H_{u_{i-1}}} \ast [0,h_j) \\
& = H_{u_{i-1}} \ast \left(  [0, h_{u_{i-1} + 1} )  
+ h_{u_{i-1} + 1} \ast  [0, h_{u_{i-1} + 2} )   \right. \\
& \qquad  \left. + h_{u_{i-1} + 1} h_{u_{i-1} + 2}\ast  [0, h_{u_{i-1} + 3} ) 
+ \cdots   \right. \\
& \qquad  \left. + h_{u_{i-1} + 1} \cdots h_{u_i -1}\ast  [0, h_{u_i} ) 
  \right) \\
& = H_{u_{i-1}} \ast \left[ 0,  \frac{H_{u_i} }{ H_{u_{i-1}}} \right) \\
& = G_{i-1} \ast \left[ 0,  g_i\right)
\end{align*}
and so $\mca = (A_i)_{i \in \N}$ is the British number system constructed 
from the integer sequence $(g_i)_{i \in \N}$.

Conversely, let $\mca = (A_i)_{i \in \N}$  be a contraction of \mcb\ constructed 
from a partition $(J_i)_{i \in \N}$ of \N\ in which some set $J_{i_0}$ 
is a not a finite interval of integers.  
Let $u = \min\left(J_{i_0}\right)$ and $w = \max\left( J_{i_0} \right)$.  
Because $J_{i_0}$ is not an interval, there is a smallest integer $v$ 
such that 
\[
u < v < w
\]
and 
$[u,v-1] \subseteq J_{i_0}$, but $v \notin  I_{j_0}$.  
Because 
\[
[u,v-1] \cup \{ w\} \subseteq J_{i_0} \subseteq [u,v-1] \cup [v+1,w]
\]
and 
\[
A_{i_0} = \sum_{j \in J_{i_0}} H_{j-1}\ast [0,h_j) 
\]
we have 
\begin{align*}
H_{u-1}\ast & [0,h_u) \cup H_{w-1}\ast [0,h_w)  \subseteq A_{i_0}   \\
& \subseteq  \sum_{ j \in [u,v-1]} H_{j-1}\ast [0,h_j) +
\sum_{ j \in [v+1,w]} H_{j-1}\ast [0,h_j) \\
& \subseteq H_{u-1} \ast \left[ 0,\frac{H_{v-1}}{H_{u-1}} \right) 
+ H_v \ast \left[ 0,\frac{H_w}{H_v} \right) 
\end{align*}
Because $h_u \geq 2$ and $h_v \geq 2$, 
it follows that 
\[
H_{u-1} \in A_{i_0}
\]
and 
\[
H_{w-1} = H_{u-1}\left( \frac{H_{w-1}}{H_{u-1}} \right)  \in A_{i_0}.
\] 

The largest multiple of $H_{u-1}$ in $H_{u-1} \ast \left[ 0,H_{v-1}/H_{u-1} \right)$
is $H_{u-1} (H_{v-1}/H_{u-1} -1)$.
The smallest positive multiple of $H_{u-1}$ in $H_v\ast\left[ 0,H_w/H_v \right)$ 
is $H_v = H_{u-1}(H_v/H_{u-1})$.
The inequality  
\[
1 \leq \frac{H_{v-1}}{H_{u-1}} - 1 < \frac{H_{v-1}}{H_{u-1}} < \frac{H_v}{H_{u-1}} 
\leq  \frac{H_{w-1}}{H_{u-1}}
\]
implies that the set $A_{i_0}$ does not contain the integer $H_{u-1}(H_{v-1}/H_{u-1})$.  
In a British number system, every set consists of consecutive multiples of its smallest positive element.  
Because the set $A_{i_0}$ lacks this property, it follows that \mca\ 
is not a British number system.  
This completes the proof.  
\end{proof}

\bt             \label{deBruijn:theorem:primes}
There is a one-to-one correspondence between sequences $ (p_i)_{i \in \N}$ 
of prime numbers and indecomposable British number systems.  
Moreover, every additive system is either indecomposable or  a contraction of an  
indecomposable system.  
\et

\begin{proof}
Let \mca\ be a British number system generated by the sequence $(g_i)_{i \in \N}$, so that 
\[
\mca = \left(G_{i-1}\ast [0,g_i) \right)_{i \in \N}.
\]
Suppose that $g_k$ is composite for some $k \in \N$.  Then $g_k = rs$, 
where $r \geq 2$ and $s \geq 2$ are integers.  
Construct the sequence $(g'_i)_{i \in \N}$ as follows:
\[
g'_i = 
\begin{cases}
g_i & \text{if $i \leq k-1  $} \\
r & \text{if $ i=k $} \\
s & \text{if $ i=k+1 $} \\
g_{i-1} & \text{if $i \geq k+2  $.} 
\end{cases}
\]
Then
\[
G'_i = \prod_{j=1}^i g'_j = 
\begin{cases}
G_i & \text{if $i \leq k-1  $} \\
r G_{k-1}& \text{if $ i=k $} \\
G_k & \text{if $ i=k+1 $} \\
G_{i-1} & \text{if $i \geq k+2  $} 
\end{cases}
\]
and
\[
\mca' = \left(G'_{i-1}\ast [0,g'_i) \right)_{i \in \N} 
\]
is the British number system generated by the sequence $(g'_i)_{i \in \N}$.  
We have 
\[
G_{i-1}\ast [0,g_i) = 
\begin{cases}
G'_{i-1}\ast [0,g'_i) & \text{if $i \leq k-1  $} \\
G'_i \ast [0,g'_{i+1}) & \text{if $i \geq k+1  $.} 
\end{cases}
\]
The identity
\[
[0,g_k) = [0,rs) = [0,r) \oplus r\ast [0,s) = [0,g'_k)+ \frac{G'_k}{G_{k-1}}\ast [0,g'_{k+1})
\]
implies that 
\[
G_{k-1}\ast [0,g_k) = G'_{k-1} \ast [0,g'_k)+ G'_k \ast [0,g'_{k+1}) 
= \sum_{i\in \{k,k+1\} } G'_{i-1} \ast [0,g'_i)  
\]
and so the British number system  \mca\ is a contraction of 
 the British number system  $\mca'$.

Conversely, if \mca\ is a contraction of a British number system 
$\mca' = \left(G'_{i-1}\ast [0,g'_i) \right)_{i \in \N}$, 
then there are a positive integer $k$ and a set $I_k$ of positive integers 
with $|I_k| \geq 2$ such that 
\[
G_{k-1}\ast [0,g_k) = \sum_{i\in I_k} G'_{i-1} \ast [0,g'_i). 
\]
Therefore,
\[
g_k = \left|  G_{k-1}\ast [0,g_k) \right| 
= \left|   \sum_{i\in I_k} G'_{i-1} \ast [0,g'_i) \right| 
= \prod_{i\in I_k} g'_i.
\]
Because $|I_k| \geq 2$ and $|g'_i| \geq 2$ for all $i\in \N$, 
it follows that the integer $g_k$ is composite.
Thus, the British number system generated by $(g_i)_{i\in \N}$ is decomposable 
if and only if $g_i$ is composite for at least one $i \in \N$.  
Equivalently, the British number system generated by $(g_i)_{i\in \N}$ 
is indecomposable if and only if $(g_i)_{i\in \N}$ is a sequence of prime numbers.  
This completes the proof.
\end{proof}

Theorem~\ref{deBruijn:theorem:primes}
has also been observed by Munagi~\cite{muna05}.

\section{Limits of additive systems}     \label{deBruijn:section:limits}

Let $\mca = (A_i)_{i\in \N_0}$ be an additive system, and let $(g_i)_{i \in [1,n]} $ be a finite 
sequence of integers with $g_i \geq 2$ for all $i \in [1,n]$.  
The \emph{dilation} of \mca\ by 
the sequence $(g_i)_{i \in [1,n]} $ is the additive system defined inductively by
\[
(g_i)_{i \in [1,n]} \ast \mca = g_1 \ast \left(   (g_i)_{i \in [2,n]} \ast \mca   \right).
\]
For $n = 1$, we have 
\begin{align*}
\mca^{(1)} & =  (g_i)_{i \in [1,1]} \ast \mca  = g_1\ast \mca \\
& = [0,g_1) \cup  (g_1\ast A_i)_{i\in \N_0} \\
& = \left( A^{(1)}_i \right)_{i\in \N_0}
\end{align*}
where
\[
A^{(1)}_1 = [0, g_1)
\]
and
\[
 A^{(1)}_i = g_1\ast A_{i-1} \text{ for $i \geq 2$.}
\]
For $n = 2$, we have 
\begin{align*}
\mca^{(2)} & =  (g_i)_{i \in [1,2]} \ast \mca  = g_1 \ast (g_2 \ast \mca)  \\
& = g_1 \ast \left( [0, g_2) \cup (g_2 \ast A_i)_{i\in \N_0} \right) \\
& =  [0, g_1) \cup \left( g_1  \ast [0, g_2)  \right)  \cup (g_1g_2\ast A_i)_{i\in \N_0}  \\
& = \left( A^{(2)}_i \right)_{i\in \N_0}
\end{align*}
where
\begin{align*}
A^{(2)}_1 & = [0, g_1) \\
A^{(2)}_2 & = g_1 \ast [0,g_2) 
\end{align*}
and
\[
A^{(2)}_i = g_1g_2\ast A_{i-2} \text{ for $i \geq 3$.}
\]
For $n = 3$, we have 
\begin{align*}
g_3 \ast \mca  
& = [0, g_3) \cup (g_3\ast A_i)_{i\in \N_0} \\
g_2\ast (g_3 \ast \mca)  & = [0, g_2) \cup g_2 \ast [0, g_3) \cup (g_2g_3\ast A_i)_{i\in \N_0} 
\end{align*}
and
\begin{align*}
\mca^{(3)} & =  (g_i)_{i \in [1,3]} \ast \mca 
 = g_1 \ast  \left( g_2\ast (g_3 \ast \mca) \right)  \\
& =  [0, g_1) \cup  \left( g_1  \ast [0, g_2) \right) \cup  \left( g_1g_2 \ast [0, g_3) \right) 
\cup (g_1g_2g_3\ast A_i)_{i\in \N_0}\\
& = \left( A_i^{(3)} \right)_{i \in \N_0}
\end{align*}
where
\begin{align*}
A_1^{(3)}  & = [0,g_1) \\
A_2^{(3)}  & = g_1 \ast [0,g_2) \\
A_3^{(3)}  & = g_1g_2 \ast  [0,g_3) \\
A_i^{(3)}  & = g_1g_2g_3A_{i-3} \qquad \text{for $i \geq 4$.}
\end{align*}

\bl         \label{deBruijn:lemma:DilateSets}
Let $(g_i)_{i =1}^n$ be a sequence of integers such that $g_i \geq 2$ for all $i$.
For every additive system $\mca = (A_i)_{i \in \N}$, 
\[
\mca^{(n)} = (g_i)_{i =1}^n \ast \mca = \left(A^{(n)}_i \right)_{i \in \N}
\]
where
\[
A^{(n)}_i = g_1g_2 \cdots g_{i-1} \ast [0,g_i) \qquad \text{for $i=1,\ldots, n$}
\]
and
\[
A^{(n)}_i = g_1g_2 \cdots g_n \ast  A_{i-n-1} \qquad \text{for $i \geq n+1$.}
\]
\el

\begin{proof}
Induction on $n$.  
\end{proof}

Let $(\mca^{(n)})_{n\in \N}$ be a sequence of additive systems.
The additive system \mca\ is the \emph{limit}
\index{limit!additive system} \index{additive system!limit} 
of the sequence $(\mca^{(n)})_{n\in \N}$ if it satisfies the following condition:
The set  $S$ belongs to \mca\ if and only if $S$  belongs to $\mca^{(n)}$ 
for all sufficiently large $n$.
We write
\[
\lim_{n\rightarrow \infty} \mca^{(n)} = \mca
\]
if \mca\ is the limit of the sequence $(\mca^{(n)})_{n\in \N}$.
The following result indicates the remarkable stability 
of a British number system.

\bt 
Let $(g_i)_{i \in \N}$ be a sequence of integers such that $g_i \geq 2$ for all $i \in \N$, 
and let \mcg\ be the British number system generated by $(g_i)_{i \in \N}$.  
Let \mca\ be an additive system and let  $\mca^{(n)} = (g_i)_{i \in [1,n]} \ast \mca$.
Then
\[
\lim_{n\rightarrow \infty} \mca^{(n)} = \mcg.
\]
\et

\begin{proof}
If $S$ is a set in \mcg, then $S = g_1g_2 \cdots g_{i-1} \ast [0,g_i)$ for some $i \in \N$.  
By Lemma~\ref{deBruijn:lemma:DilateSets}, $S$ is a set in $\mca^{(n)}$ for all $n \geq i$, 
and so $S \in \lim_{n\rightarrow \infty} \mca^{(n)}$.
  
Conversely, let $S$ be a set that is in $\mca^{(n)}$ for all sufficiently large $n$.
If $S$ is finite, then $\max(S) < g_1g_2\cdots g_k$ for some integer $k$.  
If $n \geq k$ and $i \geq n+1$, then 
\[
\max \left(A^{(n)}_i\right) \geq g_1\dots g_n \geq  g_1\dots g_k 
\]
and so $S \neq A^{(n)}_i$.  Therefore, $S = A^{(n)}_i$  for some $i \leq n$, 
and so $S = g_1g_2 \cdots g_{i-1} \ast [0,g_i)$  for some $i \leq n$.  

If $T$ is an infinite set in $\mca^{(n)}$, then $T =  g_1g_2 \cdots g_n \ast  A_{i-n-1}$ 
for some $i \geq n+1$, and so $\min(T \setminus \{0\}) \geq g_1g_2 \cdots g_n \geq 2^n$.  
If  $T \in \mca^{(n)}$ for all $n \geq N$, then $\min(T \setminus \{0\}) \geq 2^n$
for all $n \geq N$,which is absurd.    It follows that the set $S$ is in $\mca^{(n)}$ 
for all sufficiently large $n$ if and only if $S$ is finite 
and $S$ is a set in the British number system 
generated by $(g_i)_{i \in \N}$.  This completes the proof.
\end{proof}

\bc
Let $(g_i)_{i \in \N}$ be a sequence of integers such that $g_i \geq 2$ for all $i \in \N$, 
and let \mcg\ be the British number system generated by $(g_i)_{i \in \N}$.  
If $\mcg_n = (g_i)_{i \in [1,n]} \ast \N_0$, then
\[
\lim_{n\rightarrow \infty} \mcg_n = \mcg.
\]
\ec

\def\cprime{$'$} \def\cprime{$'$} \def\cprime{$'$} \def\cprime{$'$}
\providecommand{\bysame}{\leavevmode\hbox to3em{\hrulefill}\thinspace}
\providecommand{\MR}{\relax\ifhmode\unskip\space\fi MR }
\providecommand{\MRhref}[2]{%
  \href{http://www.ams.org/mathscinet-getitem?mr=#1}{#2}
}
\providecommand{\href}[2]{#2}

\end{document}